\documentclass[11pt, oneside]{article}
\usepackage{amsmath,amsfonts,amssymb,amscd}
\usepackage{tikz}
\usepackage{graphicx}
\usepackage{extarrows}
\usepackage{amsthm}
\usepackage{eucal}
\usepackage{enumerate}
\usepackage{appendix}
\date{}

\begin{document}

\bibliographystyle{plain}

\newtheorem{theorem}{Theorem}[section]
\newtheorem{proposition}{Proposition}[section]
\newtheorem{corollary}{Corollary}[theorem]
\newtheorem{lemma}[theorem]{Lemma}

\title{Examples of open manifolds with almost quadratic volume growth and infinite Betti numbers}

\author{Huihong Jiang \thanks{Partially supported by Shanghai Sailing Program (No.22YF1431400) and NSF of China (No.12201411)}}

\maketitle

\begin{abstract}
  We construct a family of examples of complete $(2+n)-$dimensional ($n\ge 2$) open manifolds with positive Ricci curvature, sectional curvature bounded from below and infinite Betti numbers $b_2,b_n$, moreover its volume growth can be arbitrarily close to quadratic volume growth. Compared with some known result of finite topology for manifolds with nonnegative Ricci curvature and lower sectional curvature bound \cite{AG}, it makes sense to ask whether complete manifolds with such curvature bounds must be of finite topological type or not provided with at most quadratic volume growth.
\end{abstract}

\section{Introduction}

An inspiring theorem about the topology of Riemannian manifolds with certain prescribed curvature assumptions was came up with by Gromov in 1980s \cite{G}: the total Betti number of complete manifolds $M^n$ (either compact or noncompact) with nonnegative sectional curvature admits a uniform upper bound only depending on the dimension $n$. Naturally, one may ask: {\it Is there any bound for the Betti numbers of complete manifolds with nonnegative Ricci curvature?}

For the first Betti number, Anderson \cite{And} proved that $b_1(M^n)\le n$ with equality only if $M^n$ is a compact manifold with flat Riemannian metric; For the codimension one Betti number, Shen and Sormani \cite{SSo} showed that either $M^n$ is a flat normal bundle over a compact totally geodesic submanifold or $M^n$ has a trivial codimension one integer homology, thus $b_{n-1}(M^n)\le1$.

However, when turning attention to the Betti numbers $b_i$ ($2\le i\le n-2$) for complete manifolds with nonnegative Ricci curvature, such kind of uniform bound was showed to be impossible without further assumptions. The earliest examples were given by J. P. Sha and D. G. Yang \cite{ShY,ShY1}, which stated that there exist complete $(m+n)-$dimensional ($m\ge2,n\ge2$) compact manifolds with positive Ricci curvature and arbitrarily large Betti numbers $b_m,b_n$. These compact examples can be extended to noncompact case \cite{SW}, which said that there exist complete $(m+n)-$dimensional ($m\ge2,n\ge2$) open manifolds with nonnegative Ricci curvature, sectional curvature bounded from below, and infinite Betti numbers $b_m,b_n$.

On the other hand, the topology of open manifolds with nonnegative Ricci curvature could be put into control with some further assumptions. Let $M^n$ be a complete $n$-dimensional noncompact manifold, $p_0\in M$ a fixed point. By $B_t(p_0)$ (resp. $S_t(p_0)$) denote the ball (resp. sphere) with the center at $p_0$ and the radius $t$. Set
\[{\text{diam}}(p_0; t)=\sup_i{\text{diam}}\left(\Sigma_i, M\setminus B_{\frac12t}(p_0)\right),\]
where $\Sigma_i$ is a connected component of $S_t(p_0)$, ${\text{diam}}\left(\Sigma_i, M\setminus B_{\frac12t}(p_0)\right)$ denotes the diameter of $\Sigma_i$ in $M\setminus B_{\frac12t}(p_0)$ when $M\setminus B_{\frac12t}(p_0)$ is considered as a metric space with the induced metric and $\Sigma_i$ as its subset.

The earliest such result was stated by Abresch and Gromoll \cite{AG}: a complete $n$-dim ($n\ge3$) manifold with nonnegative Ricci curvature and sectional curvature bounded from below must be of finite topological type provided with diameter growth $o(t^{\frac1n})$. An open manifolds is said to be of finite topological type if it is homeomorphic to the interior of a compact manifold with boundary; otherwise, infinite topological type. Obviously, the Betti numbers for manifolds of finite topological type must be finite. Later on, J. P. Sha and Z. M. Shen \cite{SS} proved that open manifolds with nonnegative Ricci curvature and slow volume growth must be of finite topological type provided with quadratically nonnegatively curved infinity. One may ask: {\it How slow does it need for the diameter growth or volume growth of complete manifolds with nonnegative Ricci curvature and sectional curvature bounded from below to be of finite topological type?}

By extending Sha-Yang's compact example \cite{ShY1} to the noncompact case, Z. M. Shen and G. F. Wei \cite{SW} gave an example with such curvature bounds of diameter growth $O(t^{\frac12+\frac1{4(m-1)/n+2}})$ and volume growth $O(t^{\frac m2+1+m\frac1{4(m-1)/n+2}})$. Topologically, their example is
\[M^{m+n}=\left((\mathbb{R}^{m+1}\setminus\coprod\limits_{i=1}^{+\infty}D_i^{m+1})\times S^{n-1}\right)
\cup_{\textrm{Id}}\coprod\limits_{i=1}^{+\infty}(S^m\times D^n)_i,\]
where by 'Id' we mean gluing along the corresponding boundaries (i.e. $S^m\times S^{n-1}$) through the identity map. Such manifold admits infinite Betti numbers $b_m,b_n$ ($m\ge2,n\ge2$). Under a similar topological construction, while with a totally different construction of the metric, we get

\begin{theorem}\label{th1}
For all integers $n\ge2$, there exists a complete $(2+n)-$dimensional Riemannian manifold $M$ with a base point $p_0$, which admits positive Ricci curvature, sectional curvature bounded from below and infinite Betti number $b_2,b_n$. Moreover, it satisfies
\[\mathrm{diam}(p_0;t)=O(t^{\frac{1+\gamma}2}), \qquad \mathrm{vol}\left(B_t(p_0)\right)=O(t^{2+\gamma}),\]
where $\gamma$ can be chosen to be a positive constant arbitrarily close to $0$.
\end{theorem}

By taking a metric product of the above example with the standard spheres, we have

\begin{corollary}\label{cor1.1}
For all integers $N\ge4$ and $2\le k\le N-2$, there exists a complete $N-$dimensional Riemannian manifold $M$ with a base point $p_0$, which admits positive Ricci curvature, sectional curvature bounded from below and infinite Betti number $b_k$. Moreover, it satisfies
\[\mathrm{diam}(p_0;t)=O(t^{\frac{1+\gamma}2}), \qquad \mathrm{vol}\left(B_t(p_0)\right)=O(t^{2+\gamma}),\]
where $\gamma$ can be chosen to be a positive constant arbitrarily close to $0$.
\end{corollary}

One may find that the constant $\gamma$ here exactly corresponds to the constant $\frac1{4(m-1)/n+2}$ in Shen-Wei's example. Both of the examples admit diameter growth arbitrarily close to degree $\frac12$ when the constant $\gamma$ (w.r.t $\frac1{4(m-1)/n+2}$) goes to $0$. The difference is that the dimension of Shen-Wei's example must go to infinity as $m\to\infty$ while the dimension of our example can be an arbitrary integer bigger than or equal to $4$ since $\gamma$ is independent of $n$.

An open manifold is said to have {\it quadratic volume growth} when $\mathrm{vol}\left(B_t(p_0)\right)=O(t^2)$. For a complete manifold with at most quadratic volume growth, S. Y. Cheng and S. T. Yau \cite{CY} proved that it does not admit any negative subharmonic function and satisfies some more properties which are concluded to be the parabolicity. On the other hand, J. Lott and Z. Shen \cite{JS} showed that manifolds with lower quadratic curvature decay and volume growth slower than quadratic volume growth must be of finite topological type when it does not collapse at infinity. Thus the geometric property of at most quadratic volume growth may give some restrictions to the topology of complete manifolds.

We shall point out that though there is a gap between the degree of diameter growth in our example (arbitrarily close to $\frac12$) with that of the result from  Abresch and Gromoll (smaller than $\frac1n$), there seems to be little gap between our example with A-G's result in view of the volume growth. More precisely, our example can admit arbitrarily close to quadratic volume growth (i.e. $O(t^{2+\gamma})$) in any dimension $N\ge4$, while the volume growth of manifolds in A-G's result must be slower than quadratic volume growth (i.e. $o(t^{\frac{n-1}n+1})$, close to quadratic volume growth as $n\to\infty$). From above, the following question may make some sense: {\it Is a complete open manifold $M^n$ ($n\ge4$) with nonnegative Ricci curvature, sectional curvature bounded from below and at most quadratic volume growth always of finite topological type ?}

An essential tool through our construction was raised by Perelman (cf. \cite{P, BWW}):
\vspace{+0.2cm}

\noindent {\bf Gluing Criterion:} {\it Let $M_{1}$, $M_{2}$ be two compact smooth manifolds of positive Ricci curvature, with isometric boundaries $\partial M_{1} \backsimeq \partial M_{2} = X$. Suppose that the normal curvatures of $\partial M_{1}$ are bigger than the negative of the normal curvatures of $\partial M_{2}$. Then, the glued manifold $M_{1}\cup_{X} M_{2}$ along the boundary $X$ can be smoothed near $X$ to produce a manifold of positive Ricci curvature.}
\vspace{+0.2cm}

So our construction will be divided into several parts, each part is a manifold with boundary. By carefully choosing some functions and constants in the metric of each part, we can make sure that each part admits positive Ricci curvature, corresponding isometric boundaries and normal curvatures satisfying the gluing criterion. According to the gluing criterion, we can glue these parts together to get the manifold we want.

\section{The construction of the metrics}

The construction of the manifold will be completed in several steps. We will give a brief introduction to each step here and present the details later.
\begin{enumerate}[(1)]
  \item We will begin with an almost warped product metric on $[t_1,\infty)\times S^2$ (denoted by $Q$). With "almost" we mean that the metric is given as
      \[ds^2=dt^2+u^2(t)\left[dx^2+f^2(t,x)d\sigma^2\right],\]
      where $d\sigma^2$ is the standard metric of the sphere $S^1$, and $f(t,x)$ is mostly $R_0\sin x$ ($R_0$ is a small number determined later) in $[0,\pi]$ with some adjustment near $x=0$ and $x=\pi$ for the smoothness of the metric. For each interval $[t_i,t_i+2]$ ($\{t_i\}_{i=1}^{\infty}$ is an increasing sequence with $t_i+2<t_{i+1}$), the metric is actually part of a spherical metric almost like
      \[dt^2+\left(\frac1{\Lambda_i}\sin\left(\Lambda_i t\right)\right)^2\left[dx^2+(R_0\sin x)^2d\sigma^2\right].\]
      Removing a geodesic ball $B_{\frac45}(o_i)$ with the center at $o_i=(t_i+1,0)$, we will get a boundary $\partial B_{\frac45}(o_i)$ which is topologically $S^2$. By carefully choosing some parameters, the resulting manifold $Q\setminus\left(\coprod\limits_{i=1}^{+\infty}B_{\frac45}(o_i)\right)$ admits positive Ricci curvature. After scaling, the boundary $\partial B_{\frac45}(o_i)$ with intrinsic metric can be verified to be a rotationally symmetric metric (denoted by $g_i$) on $S^2$ with sectional curvature bigger than $1$ and strong concavity;
  \item Though there is no explicit expression for the metric $g_i$, a technical construction given by Perelman \cite{P} allows us to glue a neck (topologically $S^2\times[0,1]$ with positive Ricci curvature and two boundary components) along the boundary $\partial B_{\frac45}(o_i)$ such that the new boundary is a round sphere with weak concavity. For each boundary component of $Q\setminus\left(\coprod\limits_{i=1}^{+\infty}B_{\frac45}(o_i)\right)$ (including the boundary $t=t_1$), we do the similar thing to get a new manifold with boundary (denoted by $Q'$);
  \item Taking a metric product of the manifold $Q'$ with round spheres $S^{n-1}$ ($n\ge2$) of radius $c_0$, each boundary component of $Q'\times S^{n-1}$ is $S^2\times S^{n-1}$ with weak concavity (denoted by $\partial D^3_i\times S^{n-1}(c_0)$). By constructing a metric with positive Ricci curvature on $S^2\times D^n$ such that its boundary is convex and isometric to $\partial D^3_i\times S^{n-1}(c_0)$ with the convexity bigger than the concavity of $\partial D^3_i\times S^{n-1}(c_0)$, we can glue some $S^2\times D^n$ together with $Q'\times S^{n-1}$ along the boundary for each $i$, thus the resulting manifold is topologically
      \[M^{2+n}=\left((\mathbb{R}^{3}\setminus\coprod\limits_{i=0}^{+\infty}D_i^{3})\times S^{n-1}\right)\cup_{\textrm{Id}}\coprod\limits_{i=0}^{+\infty}(S^2\times D^n)_i,\]
      which admits infinite Betti numbers $b_2,b_n$.
\end{enumerate}

\subsection{The construction of $Q$}

We equip $Q=[t_1,+\infty)\times S^2$ with the metric
\[ds^2=dt^2+u^2(t)\left[dx^2+f^2(t,x)d\sigma^2\right],\]
where $d\sigma^2$ is the standard metric of the sphere $S^1$. Let $T, X, \Sigma$ be an orthonormal basis of the tangent space corresponding to the directions $dt, dx, d\sigma$ respectively.

Given some constants $c_1,c_2$ and $\gamma$ satisfying
\[0<c_2<c_1<1,\qquad 0<\gamma<\frac12.\]
where $c_2$ can be arbitrarily close to $c_1$ and $\gamma$ can be arbitrarily close to $0$. Then we can define a sequence of positive numbers $\{t_i\}_{i=1}^{+\infty}$ as
\[t_{i+1}=\left[\frac4{(1-\gamma)c_2^2}t_i^{1-\gamma}\right]^{\frac1\gamma}t_i.\]
Provided with $t_1$ sufficiently large (of course bigger than $1$), we can always find such an increasing sequence $\{t_i\}_{i=1}^{+\infty}$ satisfying $t_i+2<t_{i+1}$.

\subsubsection{Construction of $u(t)$}

The function $u(t)$ is also defined by induction together with some inductive sequences $\{\Lambda_i\}_{i=1}^{+\infty}$ and $\{\psi_i\}_{i=1}^{+\infty}$.

Denote $u_i=u(t_i), u'_i=u_t(t_i)$. Let
\[u(t_1)=c_1 t_1^{\frac{1+\gamma}2},\qquad u_{t}(t_1)=\frac{c_1}{t_{1}^{\frac{1-\gamma}2}},\]

let
\[u(t)=\frac1{\Lambda_i}\sin\left[\Lambda_i(t-t_i+\psi_i)\right], \qquad t_i<t<t_i+2,\]
and
\[u_t(t)=u_t(t_i+2)\left(\frac{t_i+2}{t}\right)^{\frac{1-\gamma}2}, \qquad t_i+2<t<t_{i+1}.\]

By choosing $\Lambda_i$ and $\psi_i$ to satisfy
\[\sin(\Lambda_i\psi_i)=\sqrt{1-(u'_i)^2},\qquad \Lambda_i=\frac{\sqrt{1-(u'_i)^2}}{u_i},\]
we can make $u(t)$ a $C^1$ function. Assuming that the following inequalities
\begin{align*}
& \Lambda_i(\psi_i+2)<\frac\pi2, &  & \Lambda_i\in\left[\frac{\sqrt{1-(u'_1)^2}}{u_i},\frac1{u_i}\right], & \\
& c_2\le \frac{u_t(t)}{t^{-\frac{1-\gamma}2}}\le c_1, & & c_2t^{\frac{1+\gamma}2}\le u(t)\le c_1t^{\frac{1+\gamma}2} &
\end{align*}
hold for $t_1\le t\le t_i$ (obviously true for $i=1$), we will show that the above inequalities also hold for $t_1\le t\le t_{i+1}$:
\begin{itemize}
  \item For $t_i< t\le t_i+2$, since
      \begin{align*}
      \frac{d}{dt}\left(\frac{u_t(t)}{t^{-\frac{1-\gamma}2}}\right) &=\frac{1-\gamma}2t^{-\frac{1+\gamma}2}\cos\left[\Lambda_i(t-t_i+\psi_i)\right]- t^{\frac{1-\gamma}2}\Lambda_i\sin\left[\Lambda_i(t-t_i+\psi_i)\right] \\
      &=t^{-\frac{1+\gamma}2}\cos\left[\Lambda_i(t-t_i+\psi_i)\right]\left\{\frac{1-\gamma}2 -t\Lambda_i\tan\left[\Lambda_i(t-t_i+\psi_i)\right]\right\} \\
      &\le t^{-\frac{1+\gamma}2}\cos\left[\Lambda_i(t-t_i+\psi_i)\right]\left\{\frac{1-\gamma}2- t_i\cdot\frac{1-(u'_i)^2}{u'_iu_i}\right\} \\
      &\le t^{-\frac{1+\gamma}2}\cos\left[\Lambda_i(t-t_i+\psi_i)\right]\left\{\frac{1-\gamma}2 -\frac{1-c_1^2t_1^{-\frac{1-\gamma}2}}{c_1^2}t_i^{1-\gamma}\right\} \\
      &<0,
      \end{align*}
      provided with $t_1\ge t_1(\gamma,c_1,c_2)$. That is to say that $\frac{u_t(t)}{t^{-\frac{1-\gamma}2}}$ is monotonic decreasing for $t_i<t\le t_i+2$, so
      \[\frac{u_t(t)}{t^{-\frac{1-\gamma}2}}\le \frac{u_t(t_i)}{t_i^{-\frac{1-\gamma}2}}\le c_1.\]

      On the other hand, write $u_t(t_i+2)$ as
      \begin{align*}
      u_t(t_i+2)&=\frac{u_t(t_i+2)}{u_t(t_i)}\frac{u_t(t_i)}{u_t(t_{i-1}+2)} \ldots\frac{u_t(t_1+2)}{u_t(t_1)}u_t(t_1), \\
      &=\frac{\cos\Lambda_i(2+\psi_i)}{\cos\Lambda_i\psi_i} \left(\frac{t_{i-1}+2}{t_i}\right)^{\frac{1-\gamma}2}\ldots \frac{\cos\Lambda_1(2+\psi_1)}{\cos\Lambda_1\psi_1}\left(\frac1{t_1}\right)^{\frac{1-\gamma}2}c_1, \\
      &=c_1\left(\frac1{t_i+2}\right)^{\frac{1-\gamma}2} \prod_{j=1}^i\frac{\frac{\cos\Lambda_j(2+\psi_j)}{\cos\Lambda_j\psi_j}} {\left(\frac{t_j}{t_j+2}\right)^{\frac{1-\gamma}2}}.
      \end{align*}
      Since for $1\le j\le i$, there holds
      \[u'_j=O(t_j^{-\frac{1-\gamma}2}),\qquad u_j=O(t_j^{\frac{1+\gamma}2}),\qquad \Lambda_j=\frac{\sqrt{1-(u'_j)^2}}{u_j}=O(t_j^{-\frac{1+\gamma}2}),\]
      thereby
      \begin{align*}
      \frac{\frac{\cos\Lambda_j(2+\psi_j)}{\cos\Lambda_j\psi_j}} {\left(\frac{t_j}{t_j+2}\right)^{\frac{1-\gamma}2}} &= \frac{\cos(2\Lambda_j)-\sin(2\Lambda_j)\tan(\Lambda_j\psi_j)}{\left(1-\frac{2}{t_j+2}\right)^{\frac{1-\gamma}2}} \\
      &\ge\frac{1-2\Lambda_j^2-2\Lambda_j\tan(\Lambda_j\psi_j)} {1-\frac{1-\gamma}{t_j+2}} =\frac{1-2\Lambda_j^2-2\Lambda^2_j\frac{u_j}{u'_j}} {1-\frac{1-\gamma}{t_j+2}} \\
      &\ge\frac{1-\frac4{u_ju'_j}} {1-\frac{1-\gamma}{t_j}} \ge\frac{1-\frac4{c_2^2t_j^\gamma}} {1-\frac{1-\gamma}{t_j}}.
      \end{align*}
      According to the choice of $\{t_i\}_{i=1}^{+\infty}$ at the beginning, we have
      \[1-\frac4{c_2^2t_{j+1}^\gamma}=1-\frac{1-\gamma}{t_j},\]
      thus
      \[\frac{u_t(t_i+2)}{(t_i+2)^{-\frac{1-\gamma}2}}=c_1\prod_{j=1}^i\frac{\frac{\cos\Lambda_j(2+\psi_j)} {\cos\Lambda_j\psi_j}} {\left(\frac{t_j}{t_j+2}\right)^{\frac{1-\gamma}2}} \ge c_1\frac{1-\frac4{c_2^2t_1^\gamma}}{1-\frac{1-\gamma}{t_i}}\ge c_2,\]
      provided with $t_1\ge t_1(\gamma,c_1,c_2)$.

      Then we have
      \[0<c_2\le\frac{u_t(t)}{t^{-\frac{1-\gamma}2}}\le c_1\]
      together with
      \[c_2t^{\frac{1+\gamma}2}\le u(t)\le c_1t^{\frac{1+\gamma}2}\]
      for $t_{i+1}<t<t_{i+1}+2$.
  \item For $t_i+2\le t\le t_{i+1}$, since
      \[\frac{u_t(t)}{t^{-\frac{1-\gamma}2}}\equiv \frac{u_t(t_i+2)}{(t_i+2)^{-\frac{1-\gamma}2}},\]
      there also holds
      \[c_2\le \frac{u_t(t)}{t^{-\frac{1-\gamma}2}}\le c_1, \qquad c_2t^{\frac{1+\gamma}2}\le u(t)\le c_1t^{\frac{1+\gamma}2}.\]
\end{itemize}

Finally, at $t=t_{i+1}$, the constants $\Lambda_{i+1}$ and $\psi_{i+1}$ are again determined by
\[\cos(\Lambda_{i+1}\psi_{i+1})=u'_{i+1},\qquad \Lambda_{i+1}=\frac{\sqrt{1-(u'_{i+1})^2}}{u_{i+1}},\]
which gives
\begin{align*}
\frac{\cos\Lambda_{i+1}(2+\psi_{i+1})}{\cos\Lambda_{i+1}\psi_{i+1}} &=\cos(2\Lambda_{i+1})-\sin(2\Lambda_{i+1})\tan(\Lambda_{i+1}\psi_{i+1}) \\
&\ge 1-2\Lambda_{i+1}^2-2\Lambda^2_{i+1}\frac{u_{i+1}}{u'_{i+1}} \\
&\ge1-\frac4{c_2^2t_{i+1}^\gamma}>0,
\end{align*}
therefor
\[\Lambda_{i+1}(\psi_{i+1}+2)<\frac\pi2.\]
So for $t_1\le t\le t_{i+1}$, we have
\[u_{tt}(t)<0,\]
thus
\[\Lambda_{i+1}=\frac{\sqrt{1-(u'_{i+1})^2}}{u_{i+1}} \in\left[\frac{\sqrt{1-(u'_1)^2}}{u_{i+1}},\frac1{u_{i+1}}\right].\]

\vspace{+0.2cm}

\noindent{\bf Conclusion}\quad There exists a $C^1$ function $u(t)$ on $[t_1,+\infty)$ satisfying
\[c_2t^{-\frac{1-\gamma}2}\le u_t(t)\le c_1t^{-\frac{1-\gamma}2},\qquad c_2t^{\frac{1+\gamma}2}\le u(t)\le c_1t^{\frac{1+\gamma}2}\]
together with some constants satisfying
\[\Lambda_i=O(\frac1{u_i})=O(t_i^{-\frac{1+\gamma}2}).\]
Especially, on each interval $[t_i,t_i+2]$,
\[u(t)=\frac1{\Lambda_i}\sin\left[\Lambda_i(t-t_i+\psi_i)\right].\]

\subsubsection{Construction of $f(t,x)$}

The function $f(t,x)$ here is actually $f_{\epsilon(t)}(x)$, which depends on the parameter $\epsilon\ll1$. Then, we will set $f_{\epsilon}(x)$ to be almost $R_0\sin x$ ($R_0$ is sufficiently small and determined when gluing in necks) and smoothed near $x=0$. The original idea here comes from \cite{P}, also see \cite{M1} and \cite{M2}.

Let $\phi\in C^2(\mathbb{R})$ satisfying
\[\phi(x)=\begin{cases}
1, & x\le0,\\
0, & x\ge1.
\end{cases}\]
Set
\[\phi_\epsilon(x)=\phi\left(\frac{x-\epsilon}{\epsilon^{\frac14}-\epsilon}\right)=\begin{cases}
1, & x\le\epsilon,\\
0, & x\ge\epsilon^{\frac14}.
\end{cases}\]
Choose $b=b(\epsilon), l=l(\epsilon)$ and $\delta=\delta(\epsilon)$ to be
\[\delta=\epsilon^2+o(\epsilon^2),\qquad lb=\sqrt{3}\epsilon+o(\epsilon),\qquad l^\epsilon=\frac{\sin lb}{lb}\frac{\epsilon}{R\sin(\epsilon+\delta)}\left(\frac{lb}{\epsilon}\right)^\epsilon,\]
which ensures that the function
\[f_{\epsilon}(x)=\begin{cases}
\frac{\sin(lx)}l, & 0\le x\le b, \vspace{+0.1cm}\\
\frac{\sin(lb)}l e^{(1-\epsilon)\log\frac xb}, & b\le x\le\epsilon, \vspace{+0.1cm}\\
R_0\sin\left(x+\delta\phi_{\epsilon}(x)\right), & \epsilon\le x\le\frac\pi2
\end{cases}\]
is $C^1$ at $x=b$ and $x=\epsilon$. Moreover, $f_{\epsilon}(x)=f_{\epsilon}(\pi-x)$ for $\frac\pi2\le x\le\pi$.

We here remark that $f_\epsilon$ on $\left[\frac\pi2,\pi\right]$ is obtained by symmetric extension at $\frac\pi2$; so, all the things done around $o_i=(t_i+1,0)$ in the following (including the surgery of gluing in necks) can and must also be done around $o'_i=(t_i+1,\pi)$ (we have to remove $B_{\frac45}(o'_i)$ since the Ricci curvatures on these geodesic balls are not necessarily positive).

Thus, we just consider the case of $0\le x\le\frac\pi2$ from now on only by requiring the geodesic balls $B_{\frac45}(o_i)$ and $B_{\frac45}(o'_i)$ disjoint, which is valid since $\frac{1}{u(t_i)}\ll1$ in the construction of \S 2.1.1.

A straightforward calculation gave

\begin{proposition}\label{prop1} (\cite{M2}, Lemma 1.22)
$\forall\eta>0$, there exists $\epsilon_0>0$ such that $\forall\epsilon<\epsilon_0, \forall x\in[0,\pi]$, one has
\[-\frac{(f_\epsilon)_{xx}}{f_\epsilon}\ge1-\eta\]
and
\[\frac{1-(f_\epsilon)_x^2}{f_\epsilon^2}\ge1-\eta.\]
\end{proposition}

which will play an important role in the positiveness of Ricci curvature. On the other hand, the following proposition is crucial for the control of the principle curvatures along the boundary $\partial B_{\frac45}(o_i)$.

\begin{proposition}\label{prop2} (\cite{M2}, Lemma 1.26)
There exists $\epsilon_0>0$ such that $\forall\epsilon<\epsilon_0, \forall x\in\left[0,\frac\pi2\right]$, one has
\[A(x)=\tan x\left|\frac{(f_{\epsilon})_x}{f_{\epsilon}}-\cot x\right|<2\epsilon.\]
\end{proposition}

Now, let $\epsilon=\epsilon(t)<\epsilon_0$ be a smooth non-increasing function satisfying
\[\epsilon(t)=\epsilon_i,\quad t_{i-1}+\frac32<t<t_i+\frac12\]
where $\{\epsilon_i\}$ is a decreasing sequence of positive constants to be determined.

Set
\[f(t,x)=f_{\epsilon(t)}(x)=
\left\{
  \begin{array}{ll}
    f_{\epsilon_i}(x), & t_{i-1}+\frac32<t<t_i+\frac12, \vspace{+0.1cm}\\
    f_{\epsilon(t)}(x), & t_i+\frac12\le t\le t_i+\frac32, \vspace{+0.1cm}\\
    f_{\epsilon_{i+1}}(x), & t_i+\frac32<t<t_{i+1}+\frac12
  \end{array}
\right.\]
and $f(t,x)=f_{\epsilon_1}(x)$ for $t_1\le t\le t_1+\frac12$. Clearly,
\[f(t,x)=R_0\sin x, \qquad \epsilon^{\frac14}\le x\le\frac\pi2.\]
To sum up, $f_t\neq0$ is possible only on $\left\{(t,x)~|~t_i+\frac12\le t\le t_i+\frac32,~0\le x\le\epsilon^{\frac14}\right\}$. Since
\begin{align*}
&d\left((t_i+1\pm\frac12,\epsilon^{\frac14}),(t_i+1,0)\right) \\
\le& d\left((t_i+1\pm\frac12,\epsilon^{\frac14}),(t_i+1\pm\frac12,0)\right) +d\left((t_i+1\pm\frac12,0),(t_i+1,0)\right) \\
\le& \int_0^{\epsilon^{\frac14}}u\left(t_i+\frac32\right)dx+\frac12\\
\le& \epsilon^{\frac14}u\left(t_i+\frac32\right)+\frac12<\frac45,
\end{align*}
provided with $\epsilon(t)$ decreasing sufficiently rapidly, we have
\[\left\{(t,x)~|~t_i+\frac12\le t\le t_i+\frac32,~0\le x\le\epsilon^{\frac14}\right\}\subset B_{\frac45}(o_i),\]
then $f_t\equiv0$ outside $B_{\frac45}(o_i)$.
\vspace{+0.2cm}

\noindent{\bf Remark}\quad We here remark that, at some discrete points, $u(t)$ and $g(t)$ are only $C^1$. However, if the manifold constructed in this manner has positive Ricci curvature on the complement of those $C^1$ parts, the manifold can then be smoothen to be a $C^2$ manifold of positive Ricci curvature; for this, one can refer to \cite{P}, also \cite{BWW, M1}.

\subsubsection{Ricci curvature and principle curvature for $Q\setminus\left(\coprod\limits_{i=1}^{+\infty}B_{\frac45}(o_i)\right)$}

Since $f_t\equiv0$ outside $B_{\frac45}(o_i)$, all the nonvanishing curvature terms for $Q\setminus\left(\coprod\limits_{i=1}^{+\infty}B_{\frac45}(o_i)\right)$ are
\begin{align*}
K(T\wedge X)&=K(T\wedge\Sigma)=-\frac{u_{tt}}u>0, \\
K(X\wedge\Sigma)&=-\frac{f_{xx}}{f u^2}-\left(\frac{u_t}u\right)^2\ge\frac{1-\eta-u_t^2}{u^2}>0
\end{align*}
from Proposition \ref{prop1} and the construction of $u(t)$ (by choosing $\eta<\frac34$). Of course, the Ricci curvature of $Q$ is always positive outside $B_{\frac45}(o_i)$.

On the other hand, the boundary $\partial B_{\frac45}(o_i)$ with induced metric satisfies Perelman's Property \cite{M1,M2}, more precisely,
\[K_{int}(X\wedge Y)>\left(\max_{Z \in S_x \partial B_{\frac45}(o_i)}\left|\mathrm{\uppercase\expandafter{\romannumeral2}}_N(Z,Z)\right|\right)^2,\]
where $N$ is the unit outward normal vector of $\partial B_{\frac45}(o_i)$, $\rm{\uppercase\expandafter{\romannumeral2}}_N$ is the corresponding second fundamental form, $K_{int}$ is the intrinsic curvature, and $S\partial B_{\frac45}(o_i)\subset T\partial B_{\frac45}(o_i)$ is the unit tangent bundle.

To verify the  Perelman's Property for the boundary $\partial B_{\frac45}(o_i)$, we rewrite $N$ as
\[N=T\cos\xi+X\sin\xi,\]
and let
\[Y=X\cos\xi-T\sin\xi.\]
Since the tangent space spanned by $T$ and $X$ is isometric to the corresponding one of $S^2\left(\frac1{\Lambda_i}\right)$, we get
\[\mathrm{\uppercase\expandafter{\romannumeral2}}_N(Y,Y)=\Lambda_i\cot(\frac45\Lambda_i),\]
the other nonvanishing second fundamental forms are
\[\mathrm{\uppercase\expandafter{\romannumeral2}}_N(\Sigma,\Sigma)=\frac{u_t}u\cos\xi+\frac{f_x}{fu}\sin\xi.\]
A standard model in the case of $f(t,x)=\sin x$ gives
\[\frac{u_t}{u}\cos\xi+\frac{\cot x}{u}\sin\xi=\Lambda_i\cot(\frac45\Lambda_i).\]
According to Proposition \ref{prop2}, we get
\begin{align*}
\left|\mathrm{\uppercase\expandafter{\romannumeral2}}_N(\Sigma,\Sigma)-\Lambda_i\cot(\frac45\Lambda_i)\right|&\le A(x)\frac{\cot x|\sin\xi|}u\\
&=A(x)\left|\Lambda_i\cot(\frac45\Lambda_i)-\frac{u_t}u\cos\xi\right|\\
&\le A(x)\left(\Lambda_i\cot(\frac45\Lambda_i)+\left|\frac{u_t}u\right|\right)\\
&\le2\epsilon_i\left(\Lambda_i\cot(\frac45\Lambda_i)+\frac{c_1}{c_2t_i}\right)\\
&\le2\epsilon_i\left(\frac54+\frac{c_1}{c_2t_i}\right)\\
&\le3\epsilon_i,
\end{align*}
provided with $t_1\ge t_1(c_1,c_2)$.

To sum up,
\[\Lambda_i\cot(\frac45\Lambda_i)\le\max_{Z \in S_x \partial B_{\frac45}(o_i)}|\mathrm{\uppercase\expandafter{\romannumeral2}}_N(Z,Z)|\le \Lambda_i\cot(\frac45\Lambda_i)+3\epsilon_i.\]

On the other hand, the sectional curvature of $\partial B_{\frac45}(o_i)$ can be calculated by the Gauss equation as
\begin{align*}
K_{int}(Y\wedge\Sigma)&=K(X\wedge\Sigma)\cos^2\xi+K(T\wedge\Sigma)\sin^2\xi +\mathrm{\uppercase\expandafter{\romannumeral2}}_N(Y,Y) \mathrm{\uppercase\expandafter{\romannumeral2}}_N(\Sigma,\Sigma)\\
&\ge\frac{1-\eta-u_t^2}{u^2}\cos^2\xi+\Lambda_i^2\sin^2\xi+ \Lambda_i\cot(\frac45\Lambda_i)\left(\Lambda_i\cot(\frac45\Lambda_i)-3\epsilon_i\right)\\
&\ge\frac{\frac34-\eta}{u_i^2}\cos^2\xi+\frac{\frac34}{u_i^2}\sin^2\xi+ \Lambda_i\cot(\frac45\Lambda_i)\left(\Lambda_i\cot(\frac45\Lambda_i)-3\epsilon_i\right)\\
&\ge\frac{\frac34-\eta-32\epsilon_i u_i^2}{u_i^2}+ \left(\Lambda_i\cot(\frac45\Lambda_i)+3\epsilon_i\right)^2\\
&>\left(\Lambda_i\cot(\frac45\Lambda_i)+3\epsilon_i\right)^2,
\end{align*}
when $\eta<\frac34$ and $\epsilon$ rapidly decreasing such that $\epsilon_i u_i^2\ll1$.

In conclusion, the boundary $\partial B_{\frac45}(o_i)$ satisfies Perelman's Property with
\[\max_{Z \in S_x \partial B_{\frac45}(o_i)} \left|\rm{\uppercase\expandafter{\romannumeral2}}_{N}(Z,Z)\right| =\Lambda_i\cot(\frac45\Lambda_i)+3\epsilon_i\triangleq h_i.\]
By choosing $t_1$ sufficiently large, we can always ensure that
\[\frac12\le h_i\le2.\]

Finally, we shall verify that the boundary $t=t_1$ also satisfies Perelman's Property. In this case, the unit outward normal vector $N$ is actually $-T$, thus the nonzero second fundamental forms are
\[\mathrm{\uppercase\expandafter{\romannumeral2}}_N(X,X)=
\mathrm{\uppercase\expandafter{\romannumeral2}}_N(\Sigma,\Sigma) =\frac{u_t(t_1)}{u(t_1)}=\frac1{t_1}.\]

On the other hand,
\[K_{int}(X\wedge\Sigma)=-\frac{f_{xx}}f\frac1{u^2(t_1)}\ge\frac{1-\eta}{c_1^2t_1^{1+\gamma}} >\frac1{t_1^2},\]
provided with $t_1\ge t_1(\gamma,c_1)$.

\subsection{Gluing along the boundary with Necks and $Q\setminus\left(\coprod\limits_{i=1}^{+\infty}B_{\frac45}(o_i)\right)$}

A neck is given by Perelman \cite{P} to construct a metric with positive Ricci curvature of $S^3\times[0,1]$ connecting the core ($\mathbb{CP}^2\setminus B(q)$ with weakly convex boundary which is a round sphere) and the ambient space whose boundary component admit strong concavity. Actually, it is a metric on $S^m\times[0,1]$,
\vspace{+0.2cm}

\noindent {\bf Neck:} {\it Let $\left(S^m,g=dt^2+B^2(t)d\sigma^2\right)$ ($m\ge 2$) be a rotationally symmetric metric, where $d\sigma^2$ is the standard metric of $S^{m-1}$. If it satisfies

{\rm(\romannumeral1)} sectional curvature $>1$;

{\rm(\romannumeral2)} $0\le t\le\pi R$ and $\max_t\{B(t)\}=r$ for some $0<r<R<1$.

Then for any $0<\rho<\frac14$ such that $r^{\frac{m-1}m}<\rho<R$, there exists a metric on $S^m\times[0,1]$ such that

1) $Ric>0$;

2) the boundary component $S^m\times\{0\}$ is concave, with normal curvatures equal to $-\lambda$, and is isometric to the round sphere $S^m(\rho\lambda^{-1})$, for some $\lambda>0$;

3) the boundary component $S^m\times\{1\}$ is strictly convex, with all its normal curvatures bigger than $1$, and is isometric to $(S^m,g)$.}
\vspace{+0.2cm}

The original version of the neck by Perelman is for $m\ge3$. But we need a neck for $m=2$ in the construction. So we explain this case here: In Perelman's construction, the only raising change from $m\ge3$ to $m=2$ is the disappearance of sectional curvature $K(\Sigma_i\wedge\Sigma_j)$ on the tangent plane corresponding to $d\sigma^2$. The possible influence from the disappearance of $K(\Sigma_i\wedge\Sigma_j)$ is breaking the positiveness of Ricci curvature (more precisely, the term $\text{Ric}(\Sigma_i,\Sigma_i)$). By checking the process of Perelman's construction, one can find that the positiveness of $\text{Ric}(\Sigma_i,\Sigma_i)$ was preserved by both $K(X\wedge\Sigma_i)$ and $K(\Sigma_i\wedge\Sigma_j)$ with lower bounds $\frac{c}{t^2}$. Thus the absence of $K(\Sigma_i\wedge\Sigma_j)$ causes no trouble with $K(X\wedge\Sigma)$ remaining unchanged. For the purpose of the completeness of our construction, we attach a detailed verification of the neck for $m=2$ to the end as Appendix \ref{appendix A}.

Thus after rescaling $\partial B_{\frac45}(o_i)$ with induced metric by  $h_i^2$ (almost $\Lambda^2_i\cot^2(\frac45\Lambda_i)$), it looks like
\[dt^2+\left(R_0\cos(\frac45\Lambda_i)\sin\frac{t}{\cos(\frac45\Lambda_i)}\right)^2d\sigma^2,\]
which owns all the characteristics of metric $g$ in the Neck with $R=\cos(\frac45\Lambda_i)$ and $r=R_0$. Then we can glue the neck (rescaled by $h_i^{-2}$) with $Q\setminus B_{\frac45}(o_i)$ along the boundary $\partial B_{\frac45}(o_i)$ such that the new boundary is actually a round sphere $S^2(\rho\lambda^{-1}h_i^{-1})$ with normal curvatures equal to $-\lambda h_i$, where $\rho$ is an arbitrary positive constant satisfying $r^{\frac12}<\rho<R$ and $\lambda$ is a positive number determined by $r$ and $\rho$. Moreover, we can glue a rescaled neck with $Q$ along the boundary $t=t_1$ since this boundary also satisfies Perelman's property.
\vspace{+0.2cm}

\noindent{\bf Conclusion}\quad Given fixed numbers $0<R_0<\frac1{20}$ and $R_0^{\frac12}<\rho<\frac14$, we can always glue $Q\setminus\left(\coprod\limits_{i=1}^{+\infty}B_{\frac45}(o_i)\right)$ with some rescaled necks along the boundaries $\partial B_{\frac45}(o_i)$ and $t=t_1$, such that the resulting manifold $Q'$ satisfying
\begin{enumerate}[1)]
  \item $Ric>0$;
  \item topologically is $\mathbb{R}^3\setminus \left(\coprod\limits_{i=0}^{+\infty}D^3_i\right)$;
  \item the boundary $\partial D^3_i$ is a round sphere $S^2(\rho\lambda^{-1}h_i^{-1})$ with normal curvatures equal to $-\lambda h_i$ ($h_i=\Lambda_i\cot(\frac45\Lambda_i)+3\epsilon_i$ and $\lambda=\lambda(r,\rho)$);
\end{enumerate}

\subsection{The construction of $(S^2\times D^n)_i$}

We equip $S^2\times D^n$ with a metric
\[d\bar{s}^2=d\bar{t}^2+\bar{u}^2(\bar{t})d\sigma^2_{S^2}+\bar{v}^2(\bar{t})d\theta^2_{S^{n-1}},\]
where $d\sigma^2_{S^2}$ and $d\theta^2_{S^{n-1}}$ are the standard metrics on $S^2$ and $S^{n-1}$, and $0\le\bar{t}\le t_0$. let $T, \{\Sigma_p\}, \{\Theta_\alpha\}, $ be an orthonormal basis of the tangent space corresponding to the directions $d\bar{t}, d\sigma^2_{S^2}, d\theta^2_{S^{n-1}}$ respectively.

Next, we will give some constants successively. Firstly, there is always a constant $0<\rho<\frac1{10}$ satisfying
\[\frac{1-(8\rho)^2}{\rho^2}\ge\frac{100n}{\tanh^2\frac{\pi}{12}}.\]
With $R_0=r=\rho^4$ in the Neck, the positive constant $\lambda$ is fixed by $\rho$, i.e.
\[\lambda=\lambda(\rho)=\lambda(n)>0.\]
Then we set
\[c_0=\frac{\tanh\frac{\pi}{12}}{16\lambda},\quad k=\frac1{2c_0},\quad t_0=\frac\pi6\cdot\frac1k,\]
which are all positive constants only depending on $n$. Finally, let
\[A_i=\frac{\rho}{\lambda h_i\cosh\frac{\pi}{12}}.\]

The functions $\bar{u}(\bar{t})$ and $\bar{v}(\bar{t})$ are given as
\[\bar{u}(\bar{t})=A_i\cosh\frac{k\bar{t}}2,\quad \bar{v}(\bar{t})=\frac{\sin(k\bar{t})}k.\]

\subsubsection{The boundary and Principle curvature}

At the boundary $\bar{t}=t_0$, the intrinsic metric is
\[\left(\rho\lambda^{-1}h_i^{-1}\right)^2d\sigma^2_{S^2}+c_0^2d\theta^2_{S^{n-1}}.\]
Let $N$ be the unit outward normal vector of $\bar{t}=t_0$ (actually $N$ is exactly $T$), the corresponding second fundamental forms are
\begin{align*}
\mathrm{\uppercase\expandafter{\romannumeral2}}_N(\Sigma_p,\Sigma_p)&=\frac{\bar{u}_{\bar{t}}}{\bar{u}}\Big{|}_{\bar{t}=t_0} =\frac{k}2\tanh\frac{\pi}{12}=4\lambda, \\
\mathrm{\uppercase\expandafter{\romannumeral2}}_N(\Theta_\alpha,\Theta_\alpha)&=\frac{\bar{v}_{\bar{t}}}{\bar{v}}\Big{|}_{\bar{t}=t_0} =\sqrt{3}k,
\end{align*}
with other terms vanishing.

\subsubsection{Sectional curvature and Ricci curvature}

As a doubly warped product metric on $S^2\times D^n$, all the nonvanishing curvature terms are
\begin{align*}
K(T,\Sigma_p,\Sigma_p,T)&=-\frac{\bar{u}_{\bar{t}\bar{t}}}{\bar{u}}=-\frac{k^2}4,\\
K(\Sigma_p,\Sigma_q,\Sigma_q,\Sigma_p)&=\frac1{\bar{u}^2}-\frac{\bar{u}_{\bar{t}}^2}{\bar{u}^2} \ge\frac{1-\bar{u}_{\bar{t}}^2}{\bar{u}^2}\Big{|}_{\bar{t}=t_0} =\frac{1-(4\rho h_i^{-1})^2}{\left(\rho\lambda^{-1}h_i^{-1}\right)^2} \ge\frac{25n}{\tanh^2\frac{\pi}{12}}\lambda^2,\\
K(T,\Theta_\alpha,\Theta_\alpha,T)&=-\frac{\bar{v}_{\bar{t}\bar{t}}}{\bar{v}}=k^2,\\
K(\Theta_\alpha,\Theta_\beta,\Theta_\beta,\Theta_\alpha)&=\frac{1}{\bar{v}^2}-\frac{\bar{v}_{\bar{t}}^2}{\bar{v}^2} =k^2,\\
K(\Sigma_p,\Theta_\alpha,\Theta_\alpha,\Sigma_p)&=-\frac{\bar{u}_{\bar{t}}}{\bar{u}}\frac{\bar{v}_{\bar{t}}}{\bar{v}} =-\frac{k^2}2\cdot\frac{\tanh\frac{k\bar{t}}2}{\tan(k\bar{t})} \ge-\frac{k^2}4.
\end{align*}
Since the constants $k,\lambda$ only depend on $n$, the sectional curvature of $(S^2\times D^n)_i$ is bounded from below.

The nonzero Ricci curvature of $(S^2\times D^n)_i$ is as follows:
\begin{align*}
\mathrm{Ric}(T,T)&=-2\frac{\bar{u}_{\bar{t}\bar{t}}}{\bar{u}}-(n-1)\frac{\bar{v}_{\bar{t}\bar{t}}}{\bar{v}} =\left(n-1-\frac12\right)k^2>0,\\
\mathrm{Ric}(\Sigma_p,\Sigma_p)&=\frac1{\bar{u}^2}-\frac{\bar{u}_{\bar{t}}^2}{\bar{u}^2}
-\frac{\bar{u}_{\bar{t}\bar{t}}}{\bar{u}}
-(n-1)\frac{\bar{u}_{\bar{t}}}{\bar{u}}\frac{\bar{v}_{\bar{t}}}{\bar{v}} \\
&\ge\frac{25n}{\tanh^2\frac{\pi}{12}}\lambda^2-\frac n4k^2 =\frac{9n}{\tanh^2\frac{\pi}{12}}\lambda^2>0,\\
\mathrm{Ric}(\Theta_\alpha,\Theta_\alpha)&=(n-2)(\frac{1}{\bar{v}^2}-\frac{\bar{v}_{\bar{t}}^2}{\bar{v}^2})
-\frac{\bar{v}_{\bar{t}\bar{t}}}{\bar{v}}
-2\frac{\bar{u}_{\bar{t}}}{\bar{u}}\frac{\bar{v}_{\bar{t}}}{\bar{v}} \ge\left(n-1-\frac12\right)k^2>0.
\end{align*}

\subsection{Gluing along the boundary with $Q'\times_{c_0}S^{n-1}$ and $(S^2\times D^n)_i$}

Taking a metric product of $Q'$ with the round sphere $S^{n-1}$ of radius $c_0$, the intrinsic metric for the boundary components of the resulting manifold $Q'\times_{c_0}S^{n-1}$ is
\[\left(\rho\lambda^{-1}h_i^{-1}\right)^2d\sigma^2_{S^2}+c_0^2d\theta^2_{S^{n-1}},\]
which is exactly isometric to the intrinsic metric for the boundary of $(S^2\times D^n)_i$. Thus we can glue $(S^2\times D^n)_i$ together with $Q'\times_{c_0}S^{n-1}$ along the boundary (the metric is $C^0$ at the boundary).

The second fundamental forms for the boundary components of $Q'\times_{c_0}S^{n-1}$ are
\begin{align*}
\mathrm{\uppercase\expandafter{\romannumeral2}}_N(\Sigma_p,\Sigma_p)&=-\lambda h_i\in\left[-2\lambda,-\frac\lambda2\right], \\
\mathrm{\uppercase\expandafter{\romannumeral2}}_N(\Theta_\alpha,\Theta_\alpha)&=0.
\end{align*}
Compared with the second fundamental forms for the boundary of $(S^2\times D^n)_i$, they satisfy the Gluing Criterion \cite{P}. So the $C^0$ metric can be smoothed to be a smooth metric with positive Ricci curvature on
\[\left(Q'\times S^{n-1}\right)\cup_{\textrm{Id}}\coprod\limits_{i=0}^{+\infty}(S^2\times D^n)_i
=\left((\mathbb{R}^{3}\setminus\coprod\limits_{i=0}^{+\infty}D_i^{3})\times S^{n-1}\right)\cup_{\textrm{Id}}\coprod\limits_{i=0}^{+\infty}(S^2\times D^n)_i.\]

\subsection{Sectional curvature, diameter growth and volume growth}

Since our example is constructed in several steps, we will investigate the sectional curvature in each steps:
\begin{enumerate}[(1)]
  \item In \S 2.1.3, we have already checked that the sectional curvature of $Q\setminus\left(\coprod\limits_{i=1}^{+\infty}B_{\frac45}(o_i)\right)$ is always positive;
  \item Then we glue the neck along $\partial B_{\frac45}(o_i)$ after rescaling by $h_i^{-2}$ with $Q\setminus\left(\coprod\limits_{i=1}^{+\infty}B_{\frac45}(o_i)\right)$ to get the manifold $Q'$, where $\frac12\le h_i\le2$. So $Q'$ also admits lower bound for sectional curvature;
  \item Obviously, the metric product of $Q'$ with the round sphere $S^{n-1}$ of radius $c_0$ has sectional curvature bounded from below;On the other hand, the sectional curvature of $(S^2\times D^n)_i$ is verified to be bounded from below in \S 2.3.2.
\end{enumerate}
In conclusion, the manifold $\left(Q'\times S^{n-1}\right)\cup_{\textrm{Id}}\coprod\limits_{i=0}^{+\infty}(S^2\times D^n)_i$ admits a lower bound for sectional curvature only depending on $n$.

Broadly speaking, our surgeries are removing geodesic balls of radius $\frac45$ with metric product of round sphere $S^{n-1}(c_0)$ and gluing with $(S^2\times D^n)_i$ which are actually in a fixed scale (since $\frac12\le h_i\le2$) on the manifold $[t_1,+\infty]\times_{u(t)}S^2\times_{c_0}S^{n-1}$, where $u(t)=O(t^\frac{1+\gamma}2)$ as $t$ goes to infinity. Thus the diameter growth and volume growth of the resulting manifold are determined by $u(t)$, i.e.
\[\mathrm{diam}(p_0;t)=O(t^{\frac{1+\gamma}2}), \qquad \mathrm{vol}\left(B_t(p_0)\right)=O(t^{2+\gamma}).\]
Since the positive constant $\gamma$ can be arbitrarily close to $0$ independent of $n$, the order of the diameter growth can be a arbitrarily close to $\frac12$, meanwhile the volume growth can be arbitrarily close to quadratic volume growth.

\appendix

\section{The construction of Neck}\label{appendix A}

For the completeness of our construction, we present the details of the neck for $m=2$ here, though it is actually a simplified version of Perelman's construction \cite{P}.
\vspace{+0.2cm}

\noindent {\bf Neck:} {\it Let $\left(S^2,g=dt^2+B^2(t)d\sigma^2\right)$ be a rotationally symmetric metric, where $d\sigma^2$ is the standard metric of $S^1$. If it satisfies

{\rm(\romannumeral1)} sectional curvature $>1$;

{\rm(\romannumeral2)} $0\le t\le\pi R$ and $\max_t\{B(t)\}=r$ for some $0<r<R<1$.

Then for any $0<\rho<\frac14$ such that $r^{\frac12}<\rho<R$, there exists a metric on $S^2\times[0,1]$ such that

1) $Ric>0$;

2) the boundary component $S^2\times\{0\}$ is concave, with normal curvatures equal to $-\lambda$, and is isometric to the round sphere $S^2(\rho\lambda^{-1})$, for some $\lambda>0$;

3) the boundary component $S^2\times\{1\}$ is strictly convex, with all its normal curvatures bigger than $1$, and is isometric to $(S^2,g)$.}
\vspace{+0.2cm}

Rewrite the initial metric
\[g=dt^2+B^2(t)d\sigma^2_{S^1}\]
as
\[g=r^2\cos^2xd\sigma^2_{S^1}+A^2(x)dx^2,\]
where $-\frac\pi2\le x\le\frac\pi2$ and $A(x)$ is a smooth positive function. Since $g$ is a complete metric on $S^2$, we have
\[\left\{
    \begin{array}{ll}
      B(0)=B(\pi R)=0, \\
      B'(0)=1,\ B'(\pi R)=-1, \\
      B''(0)=B''(\pi R)=0.
    \end{array}
  \right.\]
Together with the condition of $\max_t\{B(t)\}=r$, the function $A(x)$ satisfies
\[A(\pm\frac\pi2)=r,\qquad A'(\pm\frac\pi2)=0.\]
Moreover, since
\[\int_{-\frac\pi2}^{\frac\pi2}A(x)dx=\int_{0}^{\pi R}dt=\pi R,\]
the maximal of $A(x)$ must be bigger than or equal to $R$. From above, $A(x)$ can be expressed as
\[A(x)=r\left[1+(a_\infty-1)\eta(x)\right],\qquad -\frac\pi2\le x\le\frac\pi2,\]
where
\[a_\infty=\max\limits_x\frac{A(x)}r\ge\frac Rr>1,\qquad \max\limits_x \eta(x)=1,\]
and
\[\eta(\pm\frac\pi2)=0,\qquad \eta'(\pm\frac\pi2)=0.\]

Now, consider a metric on $S^2\times[t_0,t_\infty]$
\[ds^2=dt^2+A^2(t,x)dx^2+B^2(t,x)d\sigma_{S^1}^2,\]
let $T, X, \Sigma$ be an orthonormal basis of the tangent space corresponding to the directions $dt, dx, d\sigma_{S^1}$ respectively. Let
\[A(t,x)=tb(t)\left[1+(a(t)-1)\eta(x)\right],\qquad B(t,x)=tb(t)\cos x.\]
The nonvanishing curvatures are computed as
\begin{align*}
K(T,X,X,T)&=-\frac{A_{tt}}{A}=-\left(\frac1t+\frac{b'}{b}+\frac{\eta a'}{1+(a-1)\eta}\right)' -\left(\frac1t+\frac{b'}{b}+\frac{\eta a'}{1+(a-1)\eta}\right)^2, \\
K(T,\Sigma,\Sigma,T)&=-\frac{B_{tt}}{B}=-\left(\frac{b'}{b}\right)' -\left(\frac{b'}{b}\right)^2 -2\frac1t\frac{b'}{b}, \\
K(X,\Sigma,\Sigma,X)&=\frac1{A^2}\left(-\frac{B_{xx}}{B}+\frac{A_x}{A}\frac{B_x}{B}\right)-\frac{A_t}{A}\frac{B_t}{B} \\ &=\frac1{t^2b^2\left[1+(a-1)\eta\right]^2}\left(1-\frac{(a-1)\eta_x\tan x}{1+(a-1)\eta}\right) \\
&\qquad -\left(\frac1t+\frac{b'}{b}\right)\left(\frac1t+\frac{b'}{b}+\frac{\eta a'}{1+(a-1)\eta}\right), \\
K(T,\Sigma,\Sigma,X)&=\frac{A_t}{A^2}\frac{B_x}{B}-\frac{B_{xt}}{AB} =-\frac{a'\eta\tan x}{tb\left[1+(a-1)\eta\right]^2}.
\end{align*}
Moreover, the first term of $K(X,\Sigma,\Sigma,X)$ is actually the intrinsic curvature of the hypersurfaces $t=\text{const}$, i.e.
\[K_{int}(X,\Sigma,\Sigma,X)=\frac1{t^2b^2\left[1+(a-1)\eta\right]^2}\left(1 -\frac{(a-1)\eta_x\tan x}{1+(a-1)\eta}\right).\]

Next, let
\[\frac{b'}{b}=
\left\{
  \begin{array}{ll}
    -\beta\frac{t-t_0}{2t_0^2\ln(2t_0)}, & t_0\le t\le2t_0, \vspace{+0.2cm}\\
    -\beta\frac{\ln(2t_0)}{t\ln^2t}, & t\ge2t_0,
  \end{array}
\right.\qquad \frac{a'}{a}=-\alpha\frac{b'}{b},\]
with
\[a(t_0)=1,\qquad b(t_0)=\rho,\]
where
\[\alpha=(1+\delta)\frac{\ln a_\infty}{\ln\rho-\ln r},\qquad \beta=(1-\varepsilon)\frac{\ln\rho-\ln r}{1+1/[4\ln(2t_0)]}.\]

Notice that $\alpha>1$ since $\delta$ is a small positive number to be determined and $a_\infty\ge\frac Rr>\frac{\rho}r$. Up to now, there are some positive constants $\varepsilon, \delta, t_0, t_\infty$ to be determined (notice that constants $r, R$ and $a_\infty$ are already determined by the initial metric $g$). Firstly, we choose $\varepsilon>0$ depending on the initial metric $g$ and constant $\rho$ such that $(\frac\rho{r})^{2\varepsilon} g$ still admits sectional curvature strictly bigger than $1$, that is
\[\left(\frac r\rho\right)^{2\varepsilon}\cdot \frac1{r^2\left[1+(a_\infty-1)\eta\right]^2}\left(1 -\frac{(a_\infty-1)\eta_x\tan x}{1+(a_\infty-1)\eta}\right)>1.\]
At the maximal point of $\eta(x)$, there is $\eta_x=0$ and $\eta=1$, thus
\[\left(\frac r\rho\right)^{2\varepsilon}\cdot \frac1{r^2a_\infty^2}>1,\]
which gives
\[\alpha=(1+\delta)\frac{\ln a_\infty}{\ln\rho-\ln r}<(1+\delta)\frac{-\varepsilon(\ln r^{1/2}-\ln r)-\ln r}{\ln r^{1/2}-\ln r}<2,\]
provided with $\delta\le\frac\varepsilon2$.

Next, $t_\infty$ is determined by $t_0$ and $\delta$ such that $a(t_\infty)=a_\infty$.

Finally, the constant $t_0$ ($t_0\ge2$) will be determined when ensuring the positiveness for Ricci curvature, meanwhile the positive number $\delta$ will be determined when verifying the properties of the boundary $S^2\times\{t_\infty\}$.
\vspace{+0.2cm}

\noindent{\bf Claim 1}\quad The given metric on $S^2\times[t_0,t_\infty]$ admits positive Ricci curvature provided with constant $t_0$ large enough.
\vspace{+0.2cm}

\noindent{\it Proof of Claim 1}\quad Firstly, we will prove that $K_{int}(X,\Sigma,\Sigma,X)\ge\frac{c}{t^2}$ for some constant $c>1$, which is actually the key point in controlling the Ricci curvature. For this purpose, let
\begin{align*}
\phi(t)&=\ln\left[t^2K_{int}(X,\Sigma,\Sigma,X)\right] \\
&=\ln\left[\frac1{b^2\left[1+(a-1)\eta\right]^2}\left(1-\frac{(a-1)\eta\tan x}{1+(a-1)\eta}\right)\right].
\end{align*}
We will prove that $\phi(t)$ achieves its minimum at $t=t_0$ or $t=t_\infty$, equivalently there is no minimum point in $(t_0,t_\infty)$ for $\phi(t)$. More precisely, there is no point $t_i$ in $(t_0,t_\infty)$ such that $\phi_t(t_1)=0$ and $\phi_t(t_1^-)<0, \phi_t(t_1^+)>0$.

\begin{itemize}
  \item Write $\phi_t$ as
      \begin{align*}
      \phi_t&=-\frac{b'}{b}\frac1{1+(a-1)\eta}\Big[2(1-\eta)+2\eta a(1-\alpha) \\
      &\qquad -\frac{\alpha a\tan x\eta_x}{1-\eta+\tan x\eta_x+a(\eta-\tan x\eta_x)}\Big],
      \end{align*}
      \begin{itemize}
        \item When $\tan x\cdot \eta_x\ge0$ and $\eta\ge0$, since in $(t_0,t_\infty)$,
            \[-\frac{b'}{b}>0,\qquad \frac1{1+(a-1)\eta}>0,\]
            and the expression in the large bracket is decreasing in $a(t)$ thereby decreasing in $t$ (notice that $\alpha>1$ and $\eta(x)\le1$), there is no minimum point in $(t_0,t_\infty)$ for $\phi(t)$;
      \end{itemize}
  \item Write $\phi_t$ as
      \[\phi_t=-\frac{b'}{b}\left\{2-\frac{\alpha a}{(1-\eta)+\eta a}\left[2\eta+\frac{\tan x\eta_x}{1-\eta+\tan x\eta_x+a(\eta-\tan x\eta_x)}\right]\right\},\]
      \begin{itemize}
        \item When $\tan x\cdot \eta_x\ge0$ and $\eta<0$, similarly the expression in the large brackets is decreasing in $a(t)$ thereby decreasing in $t$, thus there is no minimum point in $(t_0,t_\infty)$ for $\phi(t)$;
        \item When $\tan x\cdot \eta_x<0$ and $\eta>0$, with a view to the monotonicity of the expression in the large bracket in $a$ again, there is no minimum point in $(t_0,t_\infty)$ for $\phi(t)$;
        \item When $\tan x\cdot \eta_x<0$ and $\eta\le0$, the expression in the square brackets is negative while $\frac{\alpha a}{(1-\eta)+\eta a}$ is positive, thus
            \[\phi_t>0, \qquad \text{for}~t_0<t<t_\infty.\]
      \end{itemize}
\end{itemize}

To sum up, $\phi(t)$ achieves its minimum at $t=t_0$ or $t=t_\infty$, that is
\[\phi(t)\ge\min\{\phi(t_0),~\phi(t_\infty)\}\ge\min\{-2\ln\rho,~ 2\ln t_\infty\}\ge2.\]
So we have
\[K_{int}(X,\Sigma,\Sigma,X)\ge\frac2{t^2},\qquad \text{for}~t_0\le t\le t_\infty.\]

When $t$ is large provided $t_0$ is sufficiently large, there are
\[\frac{b'}{b}=O\left(\frac1{t\ln t}\right),\qquad \frac{a'}{a}=-\alpha\frac{b'}{b}=O\left(\frac1{t\ln t}\right).\]
Moreover, a computation gives
\[\left|\left(\frac{b'}{b}\right)'\right|\le O\left(\frac1{t^2\ln t}\right),\qquad \left|\left(\frac{a'}{a}\right)'\right|\le O\left(\frac1{t^2\ln t}\right).\]

Thus
\begin{align*}
K(X,\Sigma,\Sigma,X)&=K_{int}(X,\Sigma,\Sigma,X)-\left(\frac1t+\frac{b'}{b}\right)\left(\frac1t+\frac{b'}{b}+\frac{\eta a'}{1+(a-1)\eta}\right) \\
&\ge\frac1{t^2}+O\left(\frac1{t^2\ln t}\right).
\end{align*}

On the other hand,
\begin{align*}
\left|K(T,X,X,T)\right|&=\Big{|}-\left(\frac{b'}{b}\right)'-\left(\frac{b'}{b}\right)^2-2\frac1t\frac{b'}{b} -2\frac{\eta a}{1+(a-1)\eta}\frac1t\frac{a'}{a} \\
&\qquad -\frac{\eta a}{1+(a-1)\eta}\left[\left(\frac{a'}{a}\right)'+\left(\frac{a'}{a}\right)^2+2\frac{b'}{b}\frac{a'}{a}\right]\Big{|}\le O\left(\frac1{t^2\ln t}\right), \\
\left|K(T,\Sigma,\Sigma,T)\right|&=\Big{|}-\left(\frac{b'}{b}\right)'-\left(\frac{b'}{b}\right)^2 -2\frac1t\frac{b'}{b}\Big{|} \le O\left(\frac1{t^2\ln t}\right).
\end{align*}

Then
\begin{align*}
\text{Ric}(X,X)&=K(X,\Sigma,\Sigma,X)+ K(T,X,X,T)\ge \frac1{t^2}+O\left(\frac1{t^2\ln t}\right), \\
\text{Ric}(\Sigma,\Sigma)&=K(X,\Sigma,\Sigma,X)+ K(T,\Sigma,\Sigma,T)\ge \frac1{t^2}+O\left(\frac1{t^2\ln t}\right),
\end{align*}
which are positive when $t_0$ is larger enough. 

Finally, since $\eta(x)\le1$ and $\alpha<2$, we have
\begin{align*}
\text{Ric}(T,T)&=K(T,X,X,T)+K(T,\Sigma,\Sigma,T) \\
&=-2\left(\frac{b'}{b}\right)'-2\left(\frac{b'}{b}\right)^2-4\frac1t\frac{b'}{b} -2\frac{\eta a}{1+(a-1)\eta}\frac1t\frac{a'}{a} \\
&\qquad -\frac{\eta a}{1+(a-1)\eta}\left[\left(\frac{a'}{a}\right)'+\left(\frac{a'}{a}\right)^2+2\frac{b'}{b}\frac{a'}{a}\right] \\
&=\left(2-\alpha\frac{\eta a}{(1-\eta)+\eta a}\right)\left[-\left(\frac{b'}{b}\right)'-\frac2t\frac{b'}{b}\right] -2\left(\frac{b'}{b}\right)^2 \\
&\qquad -\frac{\eta a}{1+(a-1)\eta}\left[\left(\frac{a'}{a}\right)^2+2\frac{b'}{b}\frac{a'}{a}\right] \\
&\ge\frac{(2-\alpha)\beta}4\frac{\ln t_0}{t^2\ln^2 t} +O\left(\frac1{t^2\ln^2 t}\right),
\end{align*}
While
\[\left|K(T,\Sigma,\Sigma,X)\right|=\left|\frac{\alpha a\eta\tan x}{b\left[1+(a-1)\eta\right]^2}\right|\cdot\left|\frac1t\frac{b'}{b}\right| \le c\text{Ric}(T,T)\ll\text{Ric}(X,X)\]
for some $c>0$, provided with $t_0$ large enough, thus the Ricci curvature is positive.

\vspace{+0.2cm}
\noindent{\bf Claim 2}\quad The given metric on $S^2\times[t_0,t_\infty]$ satisfies
\[\left\{
    \begin{array}{ll}
      a(t_0)=1,\ a'(t_0)=0, \\
      b(t_0)=\rho,\ b'(t_0)=0,
    \end{array}
  \right.\qquad\text{and}\qquad
  \left\{
    \begin{array}{ll}
      a(t_\infty)=a_\infty, \\
      b(t_\infty)>r.
    \end{array}
  \right.\]
After some rescaling of the metric, the boundary component $S^2\times\{0\}$ is concave with normal curvatures equal to $-\lambda$ and is isometric to the round sphere $S^2(\rho\lambda^{-1})$ for some $\lambda>0$, while the boundary component $S^2\times\{1\}$ is strictly convex with all its normal curvatures bigger than $1$ and is isometric to $(S^2,g)$, provided with constant $\delta>0$ small enough.
\vspace{+0.2cm}

\noindent{\it Proof of Claim 2}\quad The values of $a(t), b(t)$ and their derivatives at $t=t_0$ are already given at the beginning of the construction of the metric. A straightforward computation gives
\[\ln a(t_\infty)=\int_{t_0}^{t_\infty}\frac{a'}{a}dt =\ln a_\infty,\]
where $t_\infty$ is given by
\[(1+\delta)\cdot\frac{1+1/[4\ln(2t_0)]-\ln(2t_0)/\ln t_\infty}{1+1/[4\ln(2t_0)]}=1.\]

On the other hand,
\[\ln b(t_\infty)=\ln\rho+\int_{t_0}^{t_\infty}\frac{b'}{b}dt =\ln r+\varepsilon\ln\frac{\rho}{r},\]
which gives
\[\frac{b(t_\infty)}r=\left(\frac{\rho}{r}\right)^\varepsilon>1.\]

Then at the boundary $t=t_0$, the intrinsic metric is
\[(t_0\rho)^2\left[dx^2+\cos^2xd\sigma_{S^1}^2\right],\]
which is actually a round sphere of radius $t_0\rho$, and the nonzero second fundamental forms are
\begin{align*}
\mathrm{\uppercase\expandafter{\romannumeral2}}_N(X,X)&=-\frac{A_t}{A}\Big{|}_{t=t_0}=-\frac1{t_0}, \\
\mathrm{\uppercase\expandafter{\romannumeral2}}_N(\Sigma,\Sigma)&=-\frac{B_t}{B}\Big{|}_{t=t_0}=-\frac1{t_0},
\end{align*}
with the the unit outward normal vector $N=-T$.

At the boundary $t=t_\infty$, the intrinsic metric is
\[\left(\frac{t_\infty b(t_\infty)}r\right)^2\left[A^2(x)dx^2+r^2\cos^2xd\sigma_{S^1}^2\right]=\left(\frac{t_\infty b(t_\infty)}r\right)^2\cdot g,\]
and the nonzero second fundamental forms are
\begin{align*}
\mathrm{\uppercase\expandafter{\romannumeral2}}_N(X,X)&=\frac{A_t}{A}\Big{|}_{t=t_\infty} =\frac1{t_\infty}+O\left(\frac1{t_\infty\ln t_\infty}\right), \\
\mathrm{\uppercase\expandafter{\romannumeral2}}_N(\Sigma,\Sigma)&=\frac{B_t}{B}\Big{|}_{t=t_\infty} =\frac1{t_\infty}+O\left(\frac1{t_\infty\ln t_\infty}\right),
\end{align*}
with the the unit outward normal vector $N=T$.

Rescaling the metric by $\left(\frac{t_\infty b(t_\infty)}r\right)^{-2}$, then
\begin{itemize}
  \item the boundary $S^2\times\{0\}$ is concave with normal curvatures equal to $-\lambda$ and is isometric to the round sphere $S^2(\rho\lambda^{-1})$ for $\lambda=\frac{t_\infty b(t_\infty)}{rt_0}>0$;
  \item the boundary component $S^2\times\{1\}$ is isometric to $(S^2,g)$ and is strictly convex with all its normal curvatures equal to
      \[\frac{b(t_\infty)}r\left[1+O\left(\frac1{\ln t_\infty}\right)\right]>1,\]
      provided with $t_\infty$ sufficiently large (i.e. $\delta$ sufficiently small), since $\frac{b(t_\infty)}r=\left(\frac{\rho}{r}\right)^\varepsilon>1$.
\end{itemize}

\vskip.5cm
\noindent Mathematics \& Science College, Shanghai Normal University, Shanghai, 200234\\
{\it E-mail}: jianghuihong@shnu.edu.cn

\end{document}